\documentclass[10pt]{extarticle}

\usepackage[english]{babel}
\usepackage{graphicx}
\usepackage{framed}
\usepackage[normalem]{ulem}
\usepackage{indentfirst,ragged2e}
\usepackage{amsmath,amsthm,amssymb,amsfonts,wasysym,verbatim,bbm}
\usepackage{mathtools}
\usepackage{esint}
\usepackage{cancel}
\usepackage[italicdiff]{physics}
\usepackage[T1]{fontenc}
\usepackage{lmodern,mathrsfs}
\usepackage[dvipsnames]{xcolor}
\usepackage{nicematrix}
\usepackage{wrapfig}
\usepackage[inline,shortlabels]{enumitem}
\setlist{topsep=2pt,itemsep=2pt,parsep=0pt,partopsep=0pt}
\usepackage[utf8]{inputenc}
\usepackage[a4paper,top=0.7in,bottom=0.7in,left=0.5in,right=0.5in]{geometry}
\usepackage{multicol}
\usepackage{diagbox}
\usepackage[most]{tcolorbox}
\usepackage{tikz,tikz-3dplot,tikz-cd,tkz-tab,tkz-euclide,tikzsymbols,pgf,pgfplots}
\pgfplotsset{compat=newest}
\usepgfplotslibrary{fillbetween}
\usepackage{subfiles}
\graphicspath{{images/}{../images/}} 
\usepackage[backend=bibtex,style=numeric]{biblatex}
\usepackage{csquotes} 
\usepackage[colorlinks,linkcolor=Periwinkle,citecolor=ForestGreen,urlcolor=RedViolet]{hyperref}
\usepackage[nameinlink]{cleveref}

\newcommand{\hide}[1]{} 

\usepackage[showisoZ=false]{datetime2} 

\DTMnewdatestyle{newdate}{

}

\newcommand*{\utcpm}[1]{
\ifboolexpe{test{\ifnumequal{#1}{0}}}{}{\ifnum#1<0\else+\fi{#1}}}
\DTMnewzonestyle{newzone}{
    
}
\AtBeginDocument{
    \DTMsetdatestyle{newdate}
    \DTMsetzonestyle{newzone}
}

\usepackage{fancyhdr}
\DeclareDocumentCommand\ip{ l m }{\braces#1{\langle}{\rangle}{#2}} 
\DeclareDocumentCommand\floor{ l m }{\braces#1{\lfloor}{\rfloor}{#2}} 
\DeclareDocumentCommand\ceil{ l m }{\braces#1{\lceil}{\rceil}{#2}} 

\docsvlist{A,B,C,D,F,G,I,J,K,M,N,Q,R,T,U,V,W,X,Y,Z}

\docsvlist{A,B,C,D,E,F,G,H,I,J,K,L,M,N,O,P,Q,R,S,T,U,V,W,X,Y,Z}

\docsvlist{a,b,c,d,e,f,g,i,j,k,l,m,n,o,p,q,r,s,t,u,v,w,x,y,z,A,B,C,D,E,F,G,H,I,J,K,L,M,N,O,P,Q,R,S,T,U,V,W,X,Y,Z}

\docsvlist{a,b,c,d,e,f,g,h,i,j,k,l,m,n,o,q,r,t,u,w,x,y,z,A,B,C,D,E,F,G,H,I,J,K,L,M,N,O,P,Q,R,S,T,U,V,W,X,Y,Z}

\docsvlist{a,b,c,d,e,f,g,h,i,j,k,l,m,n,o,p,q,r,s,t,u,v,w,x,y,z,A,B,C,D,E,F,G,H,I,J,K,L,M,N,O,P,Q,R,S,T,U,V,W,X,Y,Z}

\newcommand{\limn}{\lim_{n\to\infty}}

\newcommand{\limsupn}{\limsup_{n\to\infty}}

\newcommand{\sumn}[1][1]{\sum_{n=#1}^\infty}

\newcommand{\emp}{\varnothing}

\newcommand{\sub}{\subseteq}

\newcommand{\cupp}{\bigcup}

\newcommand{\dt}{\,dt}

\newcommand{\E}{\opbraces{\operatorname{\mathbb{E}}}}

\newcommand*{\Radius}{\textbf{\textsf{R}}}

\newcommand*{\frontone}[1]{\textcolor{red!80!black}{#1}}
\newcommand*{\fronttwo}[1]{\textcolor{Green}{#1}}
\newcommand*{\frontend}[1]{\textcolor{blue!80!black}{#1}}

\newsavebox{\imaginarybox}
\newcommand{\backwarditalic}[1]{%
    \mbox{%
        \sbox{\imaginarybox}{#1}%
        \hskip\wd\imaginarybox
        \pdfsave
        \pdfsetmatrix{1 0 -0.2 1}%
        \llap{\usebox{\imaginarybox}}%
        \pdfrestore
    }
}

\newtheoremstyle{mystyle}{}{}{}{}{\sffamily\bfseries}{.}{ }{}
\newtheoremstyle{cstyle}{}{}{}{}{\sffamily\bfseries}{.}{ }{\thmnote{#3}}
\makeatletter
\renewenvironment{proof}[1][\proofname] {\par\pushQED{\qed}{\normalfont\sffamily\bfseries\topsep6\p@\@plus6\p@\relax #1\@addpunct{.} }}{\popQED\endtrivlist\@endpefalse}
\makeatother
\theoremstyle{mystyle}{\newtheorem{definition}{Definition}[section]}
\theoremstyle{mystyle}{\newtheorem{proposition}[definition]{Proposition}}
\theoremstyle{mystyle}{\newtheorem{theorem}[definition]{Theorem}}
\theoremstyle{mystyle}{\newtheorem{lemma}[definition]{Lemma}}
\theoremstyle{mystyle}{\newtheorem{corollary}[definition]{Corollary}}
\theoremstyle{mystyle}{\newtheorem*{remark}{Remark}}
\theoremstyle{mystyle}{\newtheorem*{remarks}{Remarks}}
\theoremstyle{mystyle}{\newtheorem*{example}{Example}}
\theoremstyle{mystyle}{}
\theoremstyle{definition}{}
\theoremstyle{cstyle}{\newtheorem*{conjecture}{}}
\newtheoremstyle{warn}{}{}{}{}{\normalfont}{}{ }{}
\theoremstyle{warn}

\newcommand{\warningsign}[1]{\tikz[scale=#1,every node/.style={transform shape}]{
\draw[-,line width={#1*0.8mm},red,fill=yellow,rounded corners={#1*2.5mm}] (0,0)--(1,{-sqrt(3)})--(-1,{-sqrt(3)})--cycle;
\node at (0,-1) {\fontsize{48}{60}\selectfont\bfseries!};}}

\makeatletter
\def\ifemptyarg#1{%
  \if\relax\detokenize{#1}\relax 
    \expandafter\@firstoftwo
  \else
    \expandafter\@secondoftwo
  \fi}
\makeatother

\tcolorboxenvironment{definition}{boxrule=0pt,boxsep=0pt,colback={red!10},left=8pt,right=8pt,enhanced jigsaw, borderline west={2pt}{0pt}{red},sharp corners,before skip=10pt,after skip=10pt,breakable}
\tcolorboxenvironment{proposition}{boxrule=0pt,boxsep=0pt,colback={Orange!10},left=8pt,right=8pt,enhanced jigsaw, borderline west={2pt}{0pt}{Orange},sharp corners,before skip=10pt,after skip=10pt,breakable}
\tcolorboxenvironment{theorem}{boxrule=0pt,boxsep=0pt,colback={blue!10},left=8pt,right=8pt,enhanced jigsaw, borderline west={2pt}{0pt}{blue},sharp corners,before skip=10pt,after skip=10pt,breakable}
\tcolorboxenvironment{lemma}{boxrule=0pt,boxsep=0pt,colback={Cyan!10},left=8pt,right=8pt,enhanced jigsaw, borderline west={2pt}{0pt}{Cyan},sharp corners,before skip=10pt,after skip=10pt,breakable}
\tcolorboxenvironment{corollary}{boxrule=0pt,boxsep=0pt,colback={violet!10},left=8pt,right=8pt,enhanced jigsaw, borderline west={2pt}{0pt}{violet},sharp corners,before skip=10pt,after skip=10pt,breakable}
\tcolorboxenvironment{proof}{boxrule=0pt,boxsep=0pt,blanker,borderline west={2pt}{0pt}{CadetBlue!80!white},left=8pt,right=8pt,sharp corners,before skip=10pt,after skip=10pt,breakable}
\tcolorboxenvironment{remark}{boxrule=0pt,boxsep=0pt,blanker,borderline west={2pt}{0pt}{Green},left=8pt,right=8pt,before skip=10pt,after skip=10pt,breakable}
\tcolorboxenvironment{remarks}{boxrule=0pt,boxsep=0pt,blanker,borderline west={2pt}{0pt}{Green},left=8pt,right=8pt,before skip=10pt,after skip=10pt,breakable}
\tcolorboxenvironment{example}{boxrule=0pt,boxsep=0pt,blanker,borderline west={2pt}{0pt}{Black},left=8pt,right=8pt,sharp corners,before skip=10pt,after skip=10pt,breakable}
\tcolorboxenvironment{examples}{boxrule=0pt,boxsep=0pt,blanker,borderline west={2pt}{0pt}{Black},left=8pt,right=8pt,sharp corners,before skip=10pt,after skip=10pt,breakable}
\tcolorboxenvironment{conjecture}{boxrule=0pt,boxsep=0pt,colback={gray!10},left=8pt,right=8pt,enhanced jigsaw, borderline west={2pt}{0pt}{gray},sharp corners,before skip=10pt,after skip=10pt,breakable}

\newenvironment{talign*}{\let\displaystyle\textstyle\csname align*\endcsname}{\endalign}

\makeatletter
\def\fsize{\dimexpr\f@size pt\relax}
\makeatother

\usepackage[explicit]{titlesec}
\titleformat{\section}{\normalfont\sffamily\Large\bfseries}{\thesection}{12pt}{#1}
\titleformat{\subsection}{\normalfont\sffamily\large\bfseries}{\thesubsection}{12pt}{#1}
\titleformat{\subsubsection}{\normalfont\sffamily\large\bfseries}{\thesubsection}{8pt}{#1}

\titlespacing*{\section}{0pt}{5pt}{5pt}
\titlespacing*{\subsection}{0pt}{5pt}{5pt}
\titlespacing*{\subsubsection}{0pt}{5pt}{5pt}

\newcommand{\Disp}{\displaystyle}

\DeclareMathAlphabet\mathbfcal{OMS}{cmsy}{b}{n}
\makeatletter
\g@addto@macro\normalsize{
\setlength\abovedisplayskip{3pt}
\setlength\belowdisplayskip{3pt}
\setlength\abovedisplayshortskip{0pt}
\setlength\belowdisplayshortskip{0pt}}
\makeatother

\makeatletter
\renewcommand\maketitle{
    \null\vspace{6mm}
    \begin{center}
        {\Huge\sffamily\bfseries\selectfont\@title}\\
            \vspace{6mm}
        {\Large\sffamily\selectfont\@author}\\
            \vspace{6mm}
        {\large\sffamily\selectfont\@date}
    \end{center}
    \vspace{6mm}
}
\makeatother

\definecolor{abstract_color}{HTML}{0A0848}
\renewenvironment{abstract}{
\begin{tcolorbox}[
    enhanced, width=6in, center, title={Abstract},
    frame hidden, colback=white, colbacktitle=white,
    fonttitle=\sffamily\bfseries, coltitle=black,
    borderline north={1.5pt}{0pt}{abstract_color},
    attach boxed title to top center={
        yshift=-0.25mm-\tcboxedtitleheight/2, yshifttext=2mm-\tcboxedtitleheight/2
    },
    boxed title style={
        boxrule=0.5mm,
        frame code={
            \path[tcb fill frame, abstract_color] ([xshift=-4mm]frame.west)--(frame.north west)--(frame.north east)--([xshift=4mm]frame.east)--(frame.south east)--(frame.south west)--cycle;
        },
        interior code={
            \path[tcb fill interior] ([xshift=-2mm]interior.west)--(interior.north west)--(interior.north east)--([xshift=2mm]interior.east)--(interior.south east)--(interior.south west)--cycle;
        }
    },
    top=\tcboxedtitleheight/2, bottom=0pt, left=0pt, right=0pt]
    }
{\end{tcolorbox}}

\title{Counting words without strictly increasing subwords of fixed length}
\author{Senan Sekhon}
\date{\DTMToday}

\begin{document}

\thispagestyle{empty}

\maketitle

\tableofcontents

\vspace{20pt}

\begin{abstract}
    In this paper, we derive exact formulas for generating functions counting the number of $n$-ary words avoiding strictly increasing subwords of length $k$, and provide some applications of these formulas.
\end{abstract}


\begin{tcolorbox}[
    colframe=black,
    title={Summary of Results},
    fonttitle=\large\sffamily\bfseries,
    center title,
]
    \subsubsection*{Terminology}
    
    Define the following:
    \begin{itemize}
        \item An \textbf{$n$-ary word} of \textbf{length} $s$ is a tuple $(a_1,a_2,\ldots,a_s)$ of $s$ elements chosen from $\{1,2,\ldots,n\}$.
        \item A \textbf{subword} of a word $(a_1,a_2,\ldots,a_s)$ is a contiguous tuple of the form $(a_i,a_{i+1},\ldots,a_j)$ where $1\le i<j\le s$.
        \item A \textbf{streak of length $k$} in a word is a strictly increasing subword of length $k$.
        \item $\Ps(n,k,s)$ is the number of $n$-ary words of length $s$ that do not contain any streak of length $k$.
        \item $f_{n,k}(z)$ is the generating function of $\Ps(n,k,s)$, i.e. $f_{n,k}(z)=\sum_{s=0}^\infty \Ps(n,k,s)z^s$.
    \end{itemize}

    \vspace{10pt}
    
    Then $f_{n,k}(z)$ is given by:
    \begin{equation*}
        \tcbhighmath[
            colback=white, colframe=NavyBlue,
        ]{
            f_{n,k}(z)
            = \frac{1}{\Disp\sum_{r=0}^n \psi_{k,r}\binom{n}{r}z^r}
            = \frac{k}{\Disp\sum_{r=1}^{k-1} \qty(1-\omega_k^{-r})\qty(1+\omega_k^r z)^n}
        }
    \end{equation*}
    
    Where:
    \begin{align*}
        \omega_k = e^{\tfrac{2\pi i}{k}}
        &&
        \psi_{k,r} = \begin{cases}
            1 & r \equiv 0 \pmod k \\
            -1 & r \equiv 1 \pmod k \\
            0 & \text{otherwise}
        \end{cases}
    \end{align*}
\end{tcolorbox}

\newpage

\section{Preliminary Definitions}

In this paper:
\begin{itemize}
    \item $\N$ denotes the set $\{1,2,3,\ldots\}$, while $\N_0$ denotes the set $\{0,1,2,3,\ldots\}$.
    \item $\binom{n}{k}$ denotes the binomial coefficient $\frac{n!}{(n-k)!k!}$, which is defined for all $n\in\N_0$ and all $k\in\Z$. If $k<0$ or $k>n$, we have $\binom{n}{k}=0$.
    \item $\omega_k$ denotes $e^{\tfrac{2\pi i}{k}}$, the $k$\textsuperscript{th} root of unity with the least positive argument.
\end{itemize}

\begin{definition}\label{psi_k_r:definition}
    Suppose $k\in\N$, $k\ge 2$ and $r\in\Z$. The symbol $\psi_{k,r}$ is given by:
    \begin{equation*}
        \psi_{k,r} = \begin{cases}
            1 & r \equiv 0 \pmod k \\
            -1 & r \equiv 1 \pmod k \\
            0 & \text{otherwise}
        \end{cases}
    \end{equation*}
\end{definition}
\begin{remark}
    If $k=2$, this reduces to $\psi_{2,r}=(-1)^r$. Other simple cases include:
    \begin{align*}
        \psi_{3,r}=\frac{2}{\sqrt{3}}\cos(\frac{2\pi r}{3}+\frac{\pi}{6})
        &&
        \psi_{4,r}=\frac{1}{\sqrt{2}}\cos(\frac{\pi r}{2}+\frac{\pi}{4})+\frac{1}{2}(-1)^r
    \end{align*}
\end{remark}

\begin{lemma}\label{psi_k_r:alt_expression}
    $\psi_{k,r}$ can be expressed as follows:
    \begin{equation*}
        \psi_{k,r} = \frac{1}{k}\sum_{s=1}^{k-1} (1-\omega_k^{-s})\omega_k^{rs}
    \end{equation*}
    In other words, $\psi_{k,r}$ is the average of the values of the function $z^r\qty(1-\frac{1}{z})$ as $z$ ranges over all the $k$\textsuperscript{th} roots of unity. The term with $s=0$ can be omitted as this is zero.
\end{lemma}
\begin{remark}
    Treating $\qty(\psi_{k,r})_{r=0}^{k-1}$ as a vector in $\C^k$, this lemma can be restated in terms of the discrete Fourier transform:
    \begin{equation*}
        (1,-1,0,\ldots,0)
        \xrightarrow{\textsf{DFT}}
        \qty(0,1-\omega_k^{-1},1-\omega_k^{-2},\ldots,1-\omega_k)
    \end{equation*}
\end{remark}
\begin{proof}
    We start by expanding the right side, including the term with $s=0$ as this is zero:
    \begin{align*}
        \frac{1}{k}\sum_{s=1}^{k-1} (1-\omega_k^{-s})\omega_k^{rs}
        &= \frac{1}{k}\sum_{s=0}^{k-1} (1-\omega_k^{-s})\omega_k^{rs} \\
        &= \frac{1}{k}\qty(\sum_{s=0}^{k-1} \omega_k^{rs} - \sum_{s=0}^{k-1} \omega_k^{(r-1)s})
    \end{align*}
    Since $\omega_k$ is a principal $k$\textsuperscript{th} root of unity, we have $\sum_{s=0}^{k-1} \omega_k^{ts}=0$ for all $t\in\{1,2,\ldots,k\}$ (and this all $t\in\Z$ that are not multiples of $k$). Thus the first sum vanishes unless $r\equiv 0 \pmod{k}$, and the second sum vanishes unless $r\equiv 1 \pmod{k}$. When they do not vanish, they reduce to $\sum_{s=0}^{k-1} 1=k$. This yields:
    \begin{align*}
        \frac{1}{k}\sum_{s=1}^{k-1} (1-\omega_k^{-s})\omega_k^{rs}
        &= \begin{cases}
            \frac{1}{k}(k-0) & r\equiv 0 \pmod{k} \\
            \frac{1}{k}(0-k) & r\equiv 1 \pmod{k} \\
            \frac{1}{k}(0-0) & \text{otherwise}
        \end{cases} \\
        &= \begin{cases}
            1 & r\equiv 0 \pmod{k} \\
            -1 & r\equiv 1 \pmod{k} \\
            0 & \text{otherwise}
        \end{cases}
    \end{align*}
    In all cases, this is equal to $\psi_{k,r}$.
\end{proof}

\begin{proposition}\label{psi_k_r:general_identity}
    Suppose $f(z)=\alpha_0+\alpha_1z+\alpha_2z^2+\cdots$ (a polynomial or formal power series), where $\alpha_0,\alpha_1,\alpha_2,\ldots\in\C$. Then we have:
    \begin{equation}\label{psi_k_r:general_identity_eq}
        \sum_{r=0}^\infty \psi_{k,r} \alpha_r z^r
        = \frac{1}{k}\sum_{s=1}^{k-1} \qty(1-\omega_k^{-s}) f\qty(\omega_k^s z)
    \end{equation}
\end{proposition}
\begin{proof}
    \begin{align*}
        \frac{1}{k}\sum_{s=1}^{k-1} \qty(1-\omega_k^{-s}) f\qty(\omega_k^s z)
        &= \frac{1}{k}\sum_{s=1}^{k-1} \qty(1-\omega_k^{-s}) \sum_{r=0}^\infty \alpha_r\qty(\omega_k^s z)^r \\
        &= \frac{1}{k}\sum_{s=1}^{k-1} \qty(1-\omega_k^{-s}) \sum_{r=0}^\infty \alpha_r\omega_k^{rs} z^r \\
        &= \sum_{r=0}^\infty \qty(\frac{1}{k}\sum_{s=1}^{k-1} \qty(1-\omega_k^{-s}) \omega_k^{rs}) \alpha_r z^r \\
        &= \sum_{r=0}^\infty \psi_{k,r} \alpha_r z^r \tag{by \Cref{psi_k_r:alt_expression}}
    \end{align*}
\end{proof}

\begin{corollary}\label{psi_k_r:identity_1}
    Suppose $n\in\N_0$ and $z\in\C$. Then:
    \begin{equation*}
        \sum_{r=0}^n \psi_{k,r}\binom{n}{r}z^r
        = \frac{1}{k}\sum_{s=1}^{k-1} \qty(1-\omega_k^{-s}) \qty(1+\omega_k^s z)^n
    \end{equation*}
    In other words, the binomial expansion of $(1+z)^n$ with the terms ``weighted'' by $\psi_{k,r}$ is equal to the average of the values of the function $\qty(1-\frac{1}{w})(1+wz)^n$ as $w$ ranges over the $k$\textsuperscript{th} roots of unity.
\end{corollary}

This follows directly by substituting $f(z)=(1+z)^n=\sum_{r=0}^n \binom{n}{r} z^r$ into \Cref{psi_k_r:general_identity}.\\

The following identity appears in \cite[3.3]{hwgould}, and a more general identity appears in \cite[\S1.2.6, eq. 25]{knuth}.

\begin{lemma}\label{binomial_identity_1}
    For all $n,r,s\in\N_0$, we have:
    \begin{equation}\label{binomial_identity_1_eq}
        \sum_{k=r}^{n-s} \binom{k}{r}\binom{n-k}{s} = \binom{n+1}{r+s+1}
    \end{equation}
\end{lemma}
\begin{remarks}\leavevmode
    \begin{itemize}
        \item The limits of the sum are not important, as long as they cover all values of $k$ where the binomial coefficients are nonzero (at least one term is zero if $k<r$ or $k>n-s$).
        \item If $r=0$ or $s=0$, this reduces to the \emph{hockey-stick identity}: $\sum_{k=r}^n \binom{k}{r}=\binom{n+1}{r+1}$.
        \item An easy way to remember this identity is with the following mnemonic:
        \begin{equation*}
            \sum \binom{\textsf{\textit{climbing up}}}{\textsf{fixed}}
            \binom{\textsf{\hspace{1pt}\backwarditalic{climbing down}\hspace{-4pt}}}{\textsf{fixed}}
            = \binom{\textsf{sum of tops}+1}{\textsf{sum of bottoms}+1}
        \end{equation*}
    \end{itemize}
\end{remarks}
\begin{proof}
    The right side of \eqref{binomial_identity_1_eq} is the number of ways to choose $r+s+1$ objects from a set of $n+1$. Suppose we arrange these $n+1$ objects in a line. We can choose our $r+s+1$ objects as follows: First, choose one object. This splits the line into two parts, one on the left with $k$ objects and the other on the right with $n-k$ objects. We then choose $r$ objects from the left and $s$ objects from the right:
    
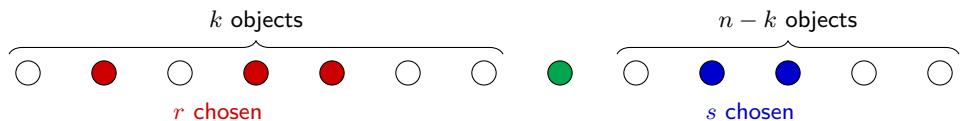
\begin{center}
        \begin{tikzpicture}
            \foreach \x in {-6,-4,-1,0,2,5,6}
                \draw (\x,0) circle (0.16);
            \foreach \x in {-5,-3,-2}
                \draw[fill=red!80!black] (\x,0) circle (0.16);
            \draw[fill=Green] (1,0) circle (0.16);
            \foreach \x in {3,4}
                \draw[fill=blue!80!black] (\x,0) circle (0.16);
            \draw[decorate,decoration={brace,amplitude=5pt}] (-6.25,0.25) -- (0.25,0.25) node[midway,above=5pt,font=\small\sffamily]{$k$ objects};
            \draw[decorate,decoration={brace,amplitude=5pt}] (1.75,0.25) -- (6.25,0.25) node[midway,above=5pt,font=\small\sffamily]{$n-k$ objects};
            \node[color=red!80!black,font=\small\sffamily] at (-3.5,-0.5) {$r$ chosen};
            \node[color=blue!80!black,font=\small\sffamily] at (3.5,-0.5) {$s$ chosen};
        \end{tikzpicture}
    \end{center}
    Summing the products yields the left side of \eqref{binomial_identity_1_eq}. Note that we overcount by a factor of $r+s+1$ as it does not matter which object we choose first (i.e. to split the line), but we overcount by the same factor on both sides of \eqref{binomial_identity_1_eq}, so these cancel out.
\end{proof}

\section{The Goulden-Jackson Cluster Method}

Suppose we have an \emph{alphabet} $\As$ and a set $\Fs\sub\cupp_{r\in\N} \As^r$ of \emph{forbidden words}. For each $r\in\N_0$, we wish to find the number of words in $\As$ of length $r$ that do not contain any of the forbidden words. For example, if $\As=\{a,b,c,d\}$ and $\Fs=\{(a,b,c),(d,a)\}$, then:
\begin{align*}
    & (a,b,d,d,c,b,d) \text{ is allowed} \\
    & (a,d,\textcolor{red!80!black}{a,b,c},b) \text{ is not allowed as it contains } (a,b,c) \\
    & (b,c,b,\textcolor{red!80!black}{d,a},b,c) \text{ is not allowed as it contains } (d,a) \\
    & (b,c,b,a,d,c) \text{ is allowed}
\end{align*}
As shown above, ``containing'' a forbidden word refers to having it as a subword, i.e. containing the letters of that word \emph{consecutively} and in the \emph{same order}.\\

We can assume that every forbidden word has length $\ge 2$, since having a single-letter forbidden word is equivalent to removing that letter from the alphabet. We can also assume that no forbidden word is a proper subword of any other forbidden word, since any word containing the longer forbidden word also contains the shorter one. Without these assumptions, the method described below still works, albeit with unnecessary complexity in the computations.\\

Suppose $\alpha_r$ is the number of words of length $r$ that do not contain any of the forbidden words. We intend to find the \textbf{generating function} $f(z)=\sum_{r=0}^\infty \alpha_r z^r$, from which we can recover the values of $\alpha_0,\alpha_1,\alpha_2,\ldots$ by $\alpha_r=\frac{f^{(r)}(0)}{r!}$.\\

This generating function can be found using the \textbf{Goulden-Jackson cluster method}, first introduced in \cite{goulden} and \cite{goulden_book}, and expanded on in \cite{noonan}. We will now proceed to describe the method.

\begin{definition}\label{front_run}
    A word $\vb{a}=(a_0,a_1,\ldots,a_{k-1})$ \textbf{front-runs} another word $\vb{b}=(b_0,b_1,\ldots,b_{l-1})$ if there exists $j\in\{0,1,\ldots,l-1\}$ such that $(a_{k-j},a_{k-j+1},\ldots,a_{k-1})=(b_0,b_1,\ldots,b_{j-1})$, i.e. the last $j$ elements of $\vb{a}$ are identical to the first $j$ elements of $\vb{b}$.\\
    The value of $j$ is known as the \textbf{overlap} of $\vb{a}$ and $\vb{b}$.
\end{definition}
\begin{example}
    $(a,b,c,b)$ front-runs $(c,b,a)$ with an overlap of $2$.
    \begin{align*}
        (a,b,\, &\textcolor{red!80!black}{c,b}) \\
        (&\textcolor{red!80!black}{c,b},a)
    \end{align*}
\end{example}
\begin{example}
    $(a,b,c,b,c)$ front-runs $(c,b,c,a,b)$ in two ways, once with an overlap of $1$ and again with an overlap of $3$.
    \begin{align*}
        (a,b,c,b,\, &\textcolor{red!80!black}{c}) & (a,b,\textcolor{red!80!black}{c,b,c}&) \\
        (&\textcolor{red!80!black}{c},b,c,a,b) & (\textcolor{red!80!black}{c,b,c}&,a,b)
    \end{align*}
\end{example}
\begin{remark}
    We do not consider a word to front-run itself in the trivial way (over its entire length). So $(a,a,a)$ front-runs itself in \emph{two} ways, with overlaps of $1$ and $2$ (and \emph{not} $3$).
\end{remark}

\begin{definition}\label{goulden_jackson:weight_equation}
    Suppose $\Fs$ is a set of forbidden words. For each forbidden word $\vb{x}\in\Fs$, the \textbf{weight} of $\vb{x}$ is given by:
    \begin{equation*}
        W(\vb{x})
        = -z^\abs{\vb{x}}
        - \sum_{\substack{\vb{y}\in\Fs \\ \vb{y} \text{ front-runs } \vb{x}}} z^{\abs{\vb{x}}-\abs{\vb{x}\wedge\vb{y}}} W(\vb{y})
    \end{equation*}
    Where $\abs{\vb{x}\wedge\vb{y}}$ is the overlap of $\vb{x}$ and $\vb{y}$. If $\vb{y}$ front-runs $\vb{x}$ in multiple ways, we must include one term for each overlap.
\end{definition}

Since the weight of a forbidden word is defined implicitly (in terms of the weights of other forbidden words), this defines a system of linear equations in the weights that needs to be solved.\\

We are now ready to state the Goulden-Jackson cluster method:

\begin{theorem}[Goulden-Jackson Cluster Method]
    Suppose $\Fs$ is a set of forbidden words, and for each $r\in\N_0$, suppose $\alpha_r$ is the number of words of length $r$ that do not contain any of the forbidden words. Then the generating function of $(\alpha_0,\alpha_1,\alpha_2,\ldots)$ is given by:
    \begin{equation*}
        f(z) = \frac{1}{1 - nz - W(\Fs)}
    \end{equation*}
    Where $W(\Fs)$ is the \textbf{weight} of $\Fs$, given by the sum of the weights of all forbidden words:
    \begin{equation*}
        W(\Fs)=\sum_{\vb{x}\in\Fs} W(\vb{x})
    \end{equation*}
\end{theorem}

\begin{example}
    Suppose $\As=\{a,b,c,d\}$ and $\Fs=\{(a,b,a),(d,a),(b,c,d)\}$. Then we have the following system of equations for the weights:
    \begin{align*}
        W(a,b,a) &= -z^3 - z^2 W(a,b,a) - z^2 W(d,a) \\
        W(d,a) &= -z^2 - z W(b,c,d) \\
        W(b,c,d) &= -z^3
    \end{align*}
    The first term on the right side of each equation is simply $-z^k$, where $k$ is the length of the forbidden word on the left. The term $-z^2W(a,b,a)$ is because $(a,b,a)$ has length $3$ and front-runs itself with an overlap of $1$, so this term gets multiplied by $z^{3-1}=z^2$. Likewise for the other terms.\\

    Solving this system, we get:
    \begin{align*}
        W(a,b,a) = \frac{-z^3 + z^4 - z^6}{1+z^2} &&
        W(d,a) = -z^2 + z^4 &&
        W(b,c,d) = -z^3
    \end{align*}
    Thus the generating function $f$ is given by:
    \begin{align*}
        f(z)
        = \frac{1}{1 - 4z - W(\Fs)}
        &= \frac{1}{1 - 4z - (W(a,b,a) + W(d,a) + W(b,c,d))} \\
        &= \frac{1}{1 - 4z - \qty(\frac{-z^3 + z^4 - z^6}{1+z^2} - z^2 + z^4 - z^3)} \\
        &= \frac{1 + z^2}{1 - 4z + 2z^2 - 2z^3 - z^4 + z^5}
    \end{align*}
    Expanding this as a power series in $z$, we get:
    \begin{equation*}
        f(z) = 1 + 4z + 15z^2 + 54z^3 + 195z^4 + 705z^5 + 2549z^6 + \cdots
    \end{equation*}
    So for example, there are $54$ words of length $3$, with letters in $\{a,b,c,d\}$, that do not contain any of $(a,b,a)$, $(d,a)$ or $(b,c,d)$. This can be verified manually: There are $4^3=64$ words in total, and we need to exclude $(a,b,a)$ and $(b,c,d)$, as well as the $4$ words that begin with $(d,a)$ and $4$ that end with $(d,a)$.
\end{example}

\begin{example}
    Suppose $\As=\{a,b,c\}$ and $\Fs=\{(a,b,c),(a,a,a)\}$. Then we have the following system of equations for the weights:
    \begin{align*}
        W(a,a,a) &= -z^3 - z^2 W(a,a,a) - z W(a,a,a) \\
        W(a,b,c) &= -z^3 - z^2 W(a,a,a)
    \end{align*}
    The last two terms in the first equation are because $(a,a,a)$ front-runs itself in two ways, with overlaps of $1$ and $2$ (the $z^2=z^{3-1}$ term comes from the overlap of $1$, and the $z=z^{3-2}$ term comes from the overlap of $2$).\\

    Solving this system, we get:
    \begin{align*}
        W(a,a,a) &= -\frac{z^3}{1+z+z^2} \\
        W(a,b,c) &= -\frac{z^3+z^4}{1+z+z^2}
    \end{align*}
    Thus the generating function $f$ is given by:
    \begin{align*}
        f(z)
        = \frac{1}{1 - 3z - W(\Fs)}
        &= \frac{1}{1 - 3z - (W(a,a,a) + W(a,b,c))} \\
        &= \frac{1}{1 - 3z - \qty(-\frac{z^3}{1+z+z^2} - \frac{z^3+z^4}{1+z+z^2})} \\
        &= \frac{1 + z + z^2}{1 - 2z - 2z^2 - z^3 + z^4}
    \end{align*}
    Expanding this as a power series in $z$, we get:
    \begin{equation*}
        f(z) = 1 + 3z + 9z^2 + 25z^3 + 70z^4 + 196z^5 + 548z^6 + \cdots
    \end{equation*}
    So for example, there are $70$ words of length $4$, with letters in $\{a,b,c\}$, that do not contain any of $(a,a,a)$ or $(a,b,c)$. To verify this, note that there are $3^4=81$ words in total, minus the $5$ words that begin or end with $(a,a,a)$ ($(a,a,a,a)$ is only counted once) and $6$ words that begin or end with $(a,b,c)$.
\end{example}

The name \emph{cluster} comes from the original terminology by Goulden and Jackson in \cite{goulden}. We avoid this terminology here by using the words \emph{front-run} and \emph{overlap} as in \Cref{front_run}.\\

In what follows, we will use the alphabet $\As=\{1,2,\ldots,n\}$ where $n\in\N$, and we will sometimes refer to the elements of $\As$ as \emph{letters} (even though they are numbers). The elements themselves do not have to be $1,2,\ldots,n$, we can use any $n$ distinct real numbers.\\

As a special case, words without the same letter appearing twice in a row are known as \emph{Smirnov words}, see \cite[\S III.24, Page 204]{flajolet} and \cite[\S1, Page 5]{iraci}.

\subsection{Radius of convergence}

In the degenerate case $\Fs=\emp$, i.e. there are no forbidden words, we have:
\begin{equation*}
    f(z) = \frac{1}{1-nz}=1+nz+n^2z^2+n^3z^3+\cdots
\end{equation*}
This makes sense as nothing is forbidden, so we are simply counting the total number of words of each length. This power series has radius of convergence $\frac{1}{n}$. We show in \Cref{radius_of_convergence:lower_bound} that in any other case ($\Fs\ne\emp$), the radius of convergence is \emph{strictly greater} than $\frac{1}{n}$.

\section{Avoiding streaks of length $k$}

In this section, we will focus on applying the Goulden-Jackson method to the case where $\Fs$ is the set of all strictly increasing words of length $k$. For example, if $n=5$ and $k=3$, we have:
\begin{align*}
    \Fs=\{(1,2,3),(1,2,4),(1,2,5),(1,3,4),(1,3,5),(1,4,5),(2,3,4),(2,3,5),(2,4,5),(3,4,5)\}
\end{align*}
In general, we have $\abs{\Fs}=\binom{n}{k}$.\\

The case $k=3$ is already proved in \cite[Theorem 3.13, Page 10]{burstein}. We will prove the general case here.

\begin{definition}
    A \textbf{streak} of length $k$ (or more simply a \emph{streak of $k$}) is any strictly increasing word of length $k$.
\end{definition}

Before proceeding to compute the weights of these forbidden words, we will prove a lemma that will greatly simplify our work.

\begin{lemma}\label{streak:w_depends_only_on_1}
    $W(a_1,a_2,\ldots,a_k)$ depends only on $a_1$ and not on $a_2,\ldots,a_k$.
\end{lemma}
\begin{proof}
    We first prove the result for $a_1=1$. Suppose $\vb{a}=(1,a_2,\ldots,a_k)$ is a streak and $\vb{b}=(b_1,b_2,\ldots,b_k)$ is a streak that front-runs $\vb{a}$. Then we have $b_i=1$ for some $2\le i\le k$. Since $\vb{b}$ is strictly increasing, this implies $b_1<1$, a contradiction. Thus there are no streaks that front-run $(1,a_2,\ldots,a_k)$, and so $W(1,a_2,\ldots,a_k)=-z^k$.\\

    Now suppose the result holds for some $a_1=p\in\{1,2,\ldots,k-1\}$. We will show that it holds for $a_1=p+1$. Suppose $\vb{a}=(p+1,a_2,\ldots,a_k)$ and $\vb{b}=(p+1,b_2,\ldots,b_k)$ are two streaks of length $k$ beginning with $p+1$. We want to show that $W(\vb{a})=W(\vb{b})$. Suppose $\vb{c}=(c_1,c_2,\ldots,c_k)$ is a streak of length $k$ that front-runs $\vb{a}$. Then we have $c_i=p+1$ for some $2\le i\le k$. Since $\vb{c}$ is strictly increasing, this implies $c_1\le p$. Now define:
    \begin{equation*}
        \vb{d} = (c_1,c_2,\ldots,c_{i-1},p+1,b_1,\ldots,b_{k-i})
    \end{equation*}
    In other words, $\vb{d}$ is a copy of $\vb{c}$ where the part overlapping with $\vb{a}$ is replaced with the corresponding elements of $\vb{b}$. Thus $\vb{d}$ front-runs $\vb{b}$ with an overlap of $k-i+1$, the same length by which $\vb{c}$ front-runs $\vb{a}$. Also, since $\vb{c}$ and $\vb{d}$ begin with the same letter $c_1\le p$, by assumption, we have $W(\vb{c})=W(\vb{d})$. Switching $\vb{a}$ and $\vb{b}$ in the above argument allows us to recover $\vb{c}$ from $\vb{d}$, with the same overlap and the same weight. Therefore, we have $W(\vb{a})=W(\vb{b})$. Thus the result holds for $p+1$, and so it holds for all $i\in\{1,2,\ldots,k\}$.
\end{proof}
\begin{remark}
    $W(a_1,a_2,\ldots,a_k)$ still depends on $k$, the length of the streak. The important part is that it does \emph{not} depend on any letters of the streak except the first. Also, this is \textbf{not} true for the general Goulden-Jackson method, but merely for the case where $\Fs$ is the set of all streaks of a fixed length.
\end{remark}

This lemma shows that we do not need to compute the weights of \emph{all} forbidden words, it suffices to use one forbidden word for each starting letter. To this end, we will fix $k\in\N$, $k\ge 2$ and define $w_m=W(m,m+1,\ldots,m+k-1)$. By \Cref{streak:w_depends_only_on_1}, this is equal to the weight of any forbidden word that starts with $m$.

\begin{example}
    Suppose $k=4$ (and $n$ remains arbitrary). We can write out the distinct weights as follows:
    \begin{align*}
        W(1,2,3,4)
        &= -z^4 \\
        W(2,3,4,5)
        &= -z^4 - zW(1,2,3,4) \\
        W(3,4,5,6)
        &= -z^4 - z^2W(1,2,3,4) - zW(1,3,4,5) - zW(2,3,4,5) \\
        W(4,5,6,7)
        &= -z^4 - z^3W(1,2,3,4) - z^2W(1,2,4,5) - z^2W(1,3,4,5) - z^2W(2,3,4,5) \\
        &\qquad\quad - zW(1,4,5,6) - zW(2,4,5,6) - zW(3,4,5,6) \\
        W(5,6,7,8)
        &= -z^4 - z^3W(1,2,3,5) - z^3W(1,2,4,5) - z^3W(1,3,4,5) - z^3W(2,3,4,5) \\
        &\qquad\quad - z^2W(1,2,5,6) - z^2W(1,3,5,6) - z^2W(1,4,5,6) \\
        &\qquad\quad - z^2W(2,3,5,6) - z^2W(2,4,5,6) - z^2W(3,4,5,6) \\
        &\qquad\quad -zW(1,5,6,7) - zW(2,5,6,7) - zW(3,5,6,7) - zW(4,5,6,7)
    \end{align*}
    Using \Cref{streak:w_depends_only_on_1}, we can simplify this to:
    \begin{align*}
        w_1
        &= -z^4 \\
        w_2
        &= -z^4 - zw_1 \\
        w_3
        &= -z^4 - z^2w_1 - zw_1 - zw_2 \\
        &= -z^4 - \qty(z^2+z)w_1 - zw_2 \\
        w_4
        &= -z^4 - z^3w_1 - z^2w_1 - z^2w_1 - z^2w_2 - zw_1 - zw_2 - zw_3 \\
        &= -z^4 - \qty(z^3+2z^2+z)w_1 - \qty(z^2+z)w_2 - zw_3 \\
        w_5
        &= -z^4 - z^3w_1 - z^3w_1 - z^3w_1 - z^3w_2 \\
        &\qquad\quad - z^2w_1 - z^2w_1 - z^2w_1 - z^2w_2 - z^2w_2 - z^2w_3 \\
        &\qquad\quad -zw_1 - zw_2 - zw_3 - zw_4 \\
        &= -z^4 - \qty(3z^3+3z^2+z)w_1 - \qty(z^3+2z^2+z)w_2 - \qty(z^2+z)w_3 - zw_4
    \end{align*}
    Continuing in this manner, we get:
    \begin{align*}
        w_6
        &= -z^4 - \qty(6z^3+4z^2+z)w_1 - \qty(3z^3+3z^2+z)w_2 - \qty(z^3+2z^2+z)w_3 - \qty(z^2+z)w_4 - zw_5 \\
        w_7
        &= -z^4 - \qty(10z^3+5z^2+z)w_1 - \qty(6z^3+4z^2+z)w_2 - \qty(3z^3+3z^2+z)w_3 \\
        &\qquad\quad - \qty(z^3+2z^2+z)w_4 - \qty(z^2+z)w_5 - zw_6
    \end{align*}
    More generally, for all $m\in\N_0$:
    \begin{talign*}
        w_m
        &= -z^4 - \qty(\binom{m-2}{2}z^3+(m-2)z^2+z)w_1 - \qty(\binom{m-3}{2}z^3+(m-3)z^2+z)w_2 - \cdots \\
        &\qquad\quad - \qty(z^3+2z^2+z)w_{m-3} - \qty(z^2+z)w_{m-2} - zw_{m-1}
    \end{talign*}
    Note that these expressions do \emph{not} depend on $n$, and so the weights will be the same for any $n\in\N$.
\end{example}

As shown above, $w_1,w_2,w_3,\ldots$ do not depend on $n$. However, the \emph{number} of weights does depend on $n$, so the generating function will also (inevitably) depend on $n$.

\begin{proposition}\label{streak:weights_recursive_formula}
    For all $m\in\N$, we have:
    \begin{equation*}
        w_m = -z^k - \sum_{r=1}^{m-1} \sum_{s=1}^{k-1} \binom{m-r-1}{s-1} z^s w_r
    \end{equation*}
\end{proposition}
\begin{proof}
    We start by expanding $w_m$ as follows:
    \begin{align*}
        w_m
        &= W(m,m+1,\ldots,m+k-1) \\
        &= -z^k \\
        &\quad -z^{k-1}\qty(W(1,2,\ldots,k-1,\frontone{m}) + \cdots + W(m-k+1,m-k+2,\ldots,\frontone{m})) \tag{$\binom{m-1}{k-1}$ terms} \\
        &\quad -z^{k-2}\qty(W(1,2,\ldots,k-2,\fronttwo{m,m+1}) + \cdots + W(m-k+2,m-k+3,\ldots,\fronttwo{m,m+1})) \tag{$\binom{m-1}{k-2}$ terms} \\
        &\quad\;\: \vdots \\
        &\quad -z^2\qty(W(1,2,\frontend{m,m+1,\ldots,m+k-3}) + \cdots + W(m-2,m-1,\frontend{m,m+1,\ldots,m+k-3})) \tag{$\binom{m-1}{2}$ terms} \\
        &\quad -z\qty(W(1,\frontend{m,m+1,\ldots,m+k-2}) + \cdots + W(m-1,\frontend{m,m+1,\ldots,m+k-2})) \tag{$m-1$ terms}
    \end{align*}
    More specifically, the terms with $z^{k-1}$ are the weights of all streaks of length $k$ that \frontone{end with $m$}. The terms with $z^{k-2}$ are the weights of all streaks of length $k$ that \fronttwo{end with $(m,m+1)$}. Likewise for the \frontend{other terms}.\\
    
    We now simplify the weights on the right using \Cref{streak:w_depends_only_on_1}. To do this, note that among the $\binom{m-1}{k-1}$ terms above with $z^{k-1}$, there are $\binom{m-2}{k-2}$ terms that start with $1$ (because they are of the form $(1,\xleftrightarrow{k-2},m)$, and the $k-2$ numbers in between are chosen from $\{2,3,\ldots,m-1\}$), $\binom{m-2}{k-3}$ terms that start with $2$, and so on. Likewise for the terms with other powers of $z$. All in all, we have:
    \begin{align*}
        w_m
        &= -z^k \\
        &\quad - z^{k-1}\qty(\binom{m-2}{k-2}w_1 + \binom{m-3}{k-2}w_2 + \cdots + \binom{k-1}{k-2}w_{m-k} + \binom{k-2}{k-2}w_{m-k+1}) \\
        &\quad - z^{k-2}\qty(\binom{m-2}{k-3}w_1 + \binom{m-3}{k-3}w_2 + \cdots + \binom{k-2}{k-3}w_{m-k+1} + \binom{k-3}{k-3}w_{m-k+2}) \\
        &\quad\;\: \vdots \\
        &\quad - z^2\qty((m-2)w_1 + (m-3)w_2 + \cdots + 2w_{m-3} + w_{m-2}) \\
        &\quad - z\qty(w_1 + w_2 + \cdots + w_{m-2} + w_{m-1})
    \end{align*}
    We now regroup the terms to separate $w_0,w_1,\ldots,w_{m-1}$:
    \begin{align*}
        w_m
        &= -z^k \\
        &\quad - \qty(\binom{m-2}{k-2}z^{k-1} + \binom{m-2}{k-3}z^{k-2} + \cdots + (m-2)z^2 + z)w_1 \\
        &\quad - \qty(\binom{m-3}{k-2}z^{k-1} + \binom{m-3}{k-3}z^{k-2} + \cdots + (m-3)z^2 + z)w_2 \\
        &\quad\;\: \vdots \\
        &\quad - \qty(z^2+z)w_{m-2} \\
        &\quad - zw_{m-1} \\
        &= -z^k - \sum_{r=1}^{m-1} \sum_{s=1}^{k-1} \binom{m-r-1}{s-1} z^s w_r
        \qedhere
    \end{align*}
\end{proof}

\begin{proposition}\label{streak:weights_general_formula}
    For all $m\in\N$, we have:
    \begin{equation}\label{streak:weights_general_formula_eq}
        w_m = -z^k\sum_{r=0}^{m-1} \psi_{k,r}\binom{m-1}{r}z^r
    \end{equation}
\end{proposition}
\begin{proof}
    When $m=1$, this reduces to $w_1=-z^k$, which we have already shown in \Cref{streak:w_depends_only_on_1}. Suppose \eqref{streak:weights_general_formula_eq} holds for $1,2,\ldots,m$. We want to show it holds for $m+1$.
    \begin{align*}
        w_{m+1}
        &= -z^k - \sum_{r=1}^m \sum_{s=1}^{k-1} \binom{m-r}{s-1} z^s w_r
            \tag{by \Cref{streak:weights_recursive_formula}} \\
        &= -z^k - \sum_{r=1}^m \sum_{s=1}^{k-1} \binom{m-r}{s-1} z^s \qty(-z^k\sum_{t=0}^{r-1} \psi_{k,t}\binom{r-1}{t}z^t)
            \tag{by assumption} \\
        &= -z^k \qty(1 - \sum_{r=1}^m \sum_{s=1}^{k-1} \binom{m-r}{s-1} z^s \qty(\sum_{t=0}^{r-1} \psi_{k,t}\binom{r-1}{t}z^t))
            \tag{factor out $-z^k$} \\
        &= -z^k \qty(1 - \sum_{r=1}^m \sum_{s=1}^{k-1} \sum_{t=0}^{r-1} \psi_{k,t}\binom{m-r}{s-1}\binom{r-1}{t}z^{s+t})
    \end{align*}
    To make things easier, we will work with $u_{m+1}=-z^{-k}w_{m+1}$ instead of $w_{m+1}$. This yields:
    \begin{align*}
        u_{m+1}
        &= 1 - \sum_{r=1}^m \sum_{s=1}^{k-1} \sum_{t=0}^{r-1} \psi_{k,t}\binom{m-r}{s-1}\binom{r-1}{t}z^{s+t} \\
        &= 1 - \sum_{r=1}^m \sum_{s=1}^{k-1} \sum_{t=0}^{r-1} \psi_{k,t}\binom{m-r}{s-1}\binom{r-1}{t}z^{s+t} \\
        &= 1 - \sum_{t=0}^{m-1} \sum_{s=1}^{k-1} \psi_{k,t} \sum_{r=t+1}^m \binom{m-r}{s-1}\binom{r-1}{t}z^{s+t}
            \tag{swap order of sums} \\
        &= 1 - \sum_{t=0}^{m-1} \sum_{s=1}^{k-1} \psi_{k,t} \binom{m}{s+t}z^{s+t}
            \tag{by \Cref{binomial_identity_1}} \\
        &= 1 - \sum_{t=0}^m \sum_{s=1}^{k-1} \psi_{k,t} \binom{m}{s+t}z^{s+t}
            \tag{the term with $t=m$ is zero} \\
        &= 1 - \underbrace{\sum_{t=0}^m \sum_{s=0}^{k-1} \psi_{k,t} \binom{m}{s+t}z^{s+t}}_{A} + \sum_{t=0}^m \psi_{k,t} \binom{m}{t}z^t
            \tag{separate the term with $s=0$}
    \end{align*}
    We now turn to the double sum $A$. We want to show that $A=1$. By the definition of $\psi_{k,t}$ (\Cref{psi_k_r:definition}), we only need to consider values of $t$ congruent to $0$ or $1 \pmod{k}$. This allows us to split the sum into two sums, setting $t=kv$ in the first and $t=kv+1$ in the second.
    \begin{equation*}
        A = \sum_{v=0}^{\floor{\frac{m}{k}}} \sum_{s=0}^{k-1} \binom{m}{s+kv}z^{s+kv} - \sum_{v=0}^{\floor{\frac{m-1}{k}}} \sum_{s=0}^{k-1} \binom{m}{s+kv+1}z^{s+kv+1}
    \end{equation*}
    Note that if $\frac{m}{k}$ is an integer, the terms in the second double sum with $v=\frac{m}{k}$ will be zero as $s+kv+1>m$. Thus we can let both sums run from $v=0$ to $\floor{\frac{m}{k}}$, which we will call $H$. This yields:
    \begin{align*}
        A
        &= \sum_{v=0}^H \sum_{s=0}^{k-1} \binom{m}{s+kv}z^{s+kv} - \sum_{v=0}^H \sum_{s=0}^{k-1} \binom{m}{s+kv+1}z^{s+kv+1} \\
        &= \sum_{v=0}^H \qty(\sum_{s=0}^{k-1} \binom{m}{s+kv}z^{s+kv} - \sum_{s=0}^{k-1} \binom{m}{s+kv+1}z^{s+kv+1}) \\
        &= \sum_{v=0}^H \qty(\sum_{s=0}^{k-1} \binom{m}{s+kv}z^{s+kv} - \sum_{s=1}^k \binom{m}{s+kv}z^{s+kv})
            \tag{$s\mapsto s-1$ in the second sum} \\
        &= \sum_{v=0}^H \qty(\binom{m}{kv}z^{kv} - \binom{m}{k+kv}z^{k+kv})
            \tag{telescoping sum} \\
        &= \sum_{v=0}^H \qty(\binom{m}{kv}z^{kv} - \binom{m}{k(v+1)}z^{k(v+1)}) \\
        &= \binom{m}{0}z^0 - \binom{m}{k(H+1)}z^{k(H+1)}
            \tag{telescoping sum} \\
        &= 1 - 0
            \tag{$k(H+1)>m$} \\
        &= 1
    \end{align*}
    Thus $A=1$. We now have:
    \begin{align*}
        u_{m+1}
        &= 1 - 1 + \sum_{t=0}^m \psi_{k,t} \binom{m}{t}z^t \\
        &= \sum_{t=0}^m \psi_{k,t} \binom{m}{t}z^t \\
        w_{m+1}
        &= -z^k\sum_{t=0}^m \psi_{k,t} \binom{m}{t}z^t
    \end{align*}
    Thus \eqref{streak:weights_general_formula_eq} holds for $m+1$, and so it holds for all $m\in\N$.
\end{proof}

\begin{proposition}\label{streak:sum_of_weights}
    The weight $W(\Fs)$ of the set of all streaks is given by:
    \begin{equation*}
        W(\Fs) = -\sum_{r=k}^n \psi_{k,r}\binom{n}{r}z^r
    \end{equation*}
\end{proposition}
\begin{proof}
    \begin{align*}
        W(\Fs)
        &= W(1,2,\ldots,k-1,k) + W(1,2,\ldots,k-1,k+1) + \cdots + W(n-k+1,n-k+2,\ldots,n-1,n) \\
        &= \binom{n-1}{k-1}w_1 + \binom{n-2}{k-1}w_2 + \cdots + \binom{k-1}{k-1}w_{n-k+1}
            \tag{by \Cref{streak:w_depends_only_on_1}} \\
        &= \sum_{s=1}^{n-k+1} \binom{n-s}{k-1}w_s \\
        &= \sum_{s=1}^{n-k+1} \binom{n-s}{k-1}\qty(-z^k\sum_{r=0}^{s-1} \psi_{k,r}\binom{s-1}{r}z^r)
            \tag{by \Cref{streak:weights_general_formula}} \\
        &= -\sum_{s=1}^{n-k+1} \sum_{r=0}^{s-1} \psi_{k,r}\binom{n-s}{k-1}\binom{s-1}{r}z^{k+r} \\
        &= -\sum_{r=0}^{n-k} \psi_{k,r}\sum_{s=r+1}^{n-k+1} \binom{n-s}{k-1}\binom{s-1}{r}z^{k+r}
            \tag{swap order of sums} \\
        &= -\sum_{r=0}^{n-k} \psi_{k,r}\binom{n}{k+r}z^{k+r}
            \tag{by \Cref{binomial_identity_1}} \\
        &= -\sum_{r=k}^n \psi_{k,r}\binom{n}{r}z^r
            \tag{$r\mapsto r-k$}
    \end{align*}
    The symbol $\psi_{k,r}$ in the last line does not change as $\psi_{k,r}=\psi_{k,r-k}$.
\end{proof}

\begin{theorem}\label{streak:generating_function}
    The generating function $f(z)$ is given by:
    \begin{equation}\label{streak:generating_function_eq}
        f(z) = \frac{1}{\Disp\sum_{r=0}^n \psi_{k,r}\binom{n}{r}z^r}
    \end{equation}
\end{theorem}
\begin{proof}
    \begin{align*}
        f(z)
        &= \frac{1}{1-nz-W(\Fs)} \\
        &= \frac{1}{\Disp 1-nz+\sum_{r=k}^n \psi_{k,r}\binom{n}{r}z^r}
            \tag{by \Cref{streak:sum_of_weights}} \\
        &= \frac{1}{\Disp\sum_{r=0}^n \psi_{k,r}\binom{n}{r}z^r}
            \tag{$\psi_{k,r}=0$ for $1<r<k$}
    \end{align*}
\end{proof}

\begin{corollary}[Alternative expression for $f(z)$]\label{streak:generating_function_alt}
    The generating function $f(z)$ is equivalently given by:
    \begin{equation}\label{streak:generating_function_alt_eq}
        f(z) = \frac{k}{\Disp\sum_{s=1}^{k-1} \qty(1-\omega_k^{-s})\qty(1+\omega_k^s z)^n}
    \end{equation}
\end{corollary}

This follows directly by substituting \Cref{psi_k_r:identity_1} into \Cref{streak:generating_function}.

\begin{remark}[Special case $k=2$]
    If $k=2$, there is a shortcut: A word does not contain a streak of length $2$ if and only if it is non-increasing. The number of non-increasing words of length $s$ is simply the number of combinations with replacement of $\{1,2,\ldots,n\}$ of length $s$, i.e. $\binom{n+s-1}{s}$. Thus the generating function is given by:
    \begin{equation*}
        f(z)
        = \sum_{s=0}^\infty \binom{n+s-1}{s}z^s
        = \frac{1}{(1-z)^n}
    \end{equation*}
    This agrees with \eqref{streak:generating_function_eq} and \eqref{streak:generating_function_alt_eq} as $\psi_{2,r}=(-1)^r$ and $\omega_2=-1$.
\end{remark}

\newcommand*{\oeis}[1]{\href{https://oeis.org/#1}{\texttt{#1}}}
The case $k=3$ appears as sequence \oeis{A225682} in OEIS \cite{oeis} as a sequence of values of $\psi_{3,r}\binom{n}{r}$. The coefficients of the power series for small values of $n$ also appear individually:

\begin{center}
    \begin{tabular}{|*{7}{c|}}
        \hline
        $n$ & 2 & 3 & 4 & 5 & 6 & 7 \\
        \hline
        \textsf{Sequence} & \oeis{A000079} & \oeis{A076264} & \oeis{A072335} & \oeis{A200781} & \oeis{A200782} & \oeis{A200783} \\
        \hline
    \end{tabular}
\end{center}

Variants of this problem involving permutations (i.e. where the letters in a word are forced to be distinct) are explored in \cite{bona}, \cite{flajolet}, \cite{kitaev} and \cite{lin}.

\section{Avoiding soft streaks of length $k$}

In this section, we state a conjecture for the case where $\Fs$ is the set of all non-decreasing words of length $k$. For example, if $n=4$ and $k=3$, we have:
\begin{align*}
    \Fs = \{&(1,1,1),(1,1,2),(1,1,3),(1,1,4),(1,2,2),(1,2,3),(1,2,4),(1,3,3),(1,3,4),(1,4,4), \\
    &(2,2,2),(2,2,3),(2,2,4),(2,3,3),(2,3,4),(2,4,4),(3,3,3),(3,3,4),(3,4,4),(4,4,4)\}
\end{align*}
In general, we have $\abs{\Fs}=\binom{n+k-1}{k}$.

\begin{definition}
    A \textbf{soft streak} of length $k$ (or more simply a \emph{soft streak of $k$}) is any non-decreasing word of length $k$.
\end{definition}

\begin{definition}
    The \textbf{generalized binomial coefficient}\footnotemark, denoted by $\Bs(n,k,r)$, is the coefficient of $x^r$ in the expansion of $\qty(1+x+x^2+\cdots+x^{k-1})^n$.
\end{definition}
\vspace*{-\fsize}
\footnotetext{This is not standard terminology, and the term ``generalized binomial coefficient'' is sometimes used with other meanings.}
\begin{remark}
    When $k=2$, this reduces to the binomial coefficient $\Bs(n,2,r)=\binom{n}{r}$.
\end{remark}

\begin{proposition}
    For all $n,k,r\in\N_0$, $k\ge 2$, we have:
    \begin{equation*}
        \Bs(n,k,r)
        = \sum_{s=0}^{\floor{\tfrac{r}{k}}} (-1)^s\binom{n}{s}\binom{n+r-ks-1}{n-1}
    \end{equation*}
\end{proposition}
\begin{proof}
    \begin{align*}
        \qty(1+x+x^2+\cdots+x^{k-1})^n
        &= \qty(1-x^k)^n(1-x)^{-n} \\
        &= \sum_{s=0}^n (-1)^s\binom{n}{s}x^{ks} \sum_{t=0}^\infty \binom{n+t-1}{n-1}x^t \\
        &= \sum_{s=0}^n (-1)^s\binom{n}{s}x^{ks} \sum_{r=ks}^\infty \binom{n+r-ks-1}{n-1}x^{r-ks}
            \tag{$r=t+ks$} \\
        &= \sum_{s=0}^n \sum_{r=ks}^\infty (-1)^s\binom{n}{s}\binom{n+r-ks-1}{n-1}x^r \\
        &= \sum_{r=0}^\infty \sum_{s=0}^{\floor{\tfrac{r}{k}}} (-1)^s\binom{n}{s}\binom{n+r-ks-1}{n-1}x^r
            \tag{swap order of sums}
    \end{align*}
    The coefficient of $x^r$ can now be read off from the inner sum.
\end{proof}

\begin{corollary}\label{psi_k_r:identity_2}
    \begin{equation*}
        \sum_{r=0}^{(k-1)n} \psi_{k,r}\Bs(n,k,r)z^r
        = \frac{\qty(1-z^k)^n}{k}\sum_{s=1}^k \qty(1-\omega_k^{-s})\qty(1-\omega_k^s z)^{-n}
    \end{equation*}
\end{corollary}
\begin{proof}
    Without the factor $\psi_{k,r}$, the left side becomes:
    \begin{equation*}
        \sum_{r=0}^{(k-1)n} \Bs(n,k,r)z^r = \qty(1+z+z^2+\cdots+z^{k-1})^n = \qty(\frac{1-z^k}{1-z})^n
    \end{equation*}
    By \Cref{psi_k_r:general_identity}, we have:
    \begin{align*}
        \sum_{r=0}^{(k-1)n} \psi_{k,r}\Bs(n,k,r)z^r
        &= \frac{1}{k}\sum_{s=1}^k \qty(1-\omega_k^{-s})\qty(\frac{1-\qty(\omega_k^s z)^k}{1-\omega_k^s z})^n \\
        &= \frac{1}{k}\sum_{s=1}^k \qty(1-\omega_k^{-s})\qty(\frac{1-z^k}{1-\omega_k^s z})^n \\
        &= \frac{\qty(1-z^k)^n}{k}\sum_{s=1}^k \qty(1-\omega_k^{-s})\qty(1-\omega_k^s z)^{-n} \qedhere
    \end{align*}
\end{proof}

\begin{conjecture}[Conjecture]\label{soft_streak:generating_function}
    The generating function $f(z)$ is given by:
    \begin{equation}\label{soft_streak:generating_function_eq}
        f(z) = \frac{\qty(1-z^k)^n}{\Disp\sum_{r=0}^{(k-1)n} \psi_{k,r}\Bs(n,k,r)z^r}
    \end{equation}
\end{conjecture}

We were not able to prove this conjecture.\\

By substituting \Cref{psi_k_r:identity_2} into \eqref{soft_streak:generating_function_eq}, we obtain an alternative expression for $f(z)$:
\begin{equation}\label{soft_streak:generating_function_alt_eq}
        f(z) = \frac{k}{\Disp\sum_{s=1}^k \qty(1-\omega_k^{-s})\qty(1-\omega_k^s z)^{-n}}
    \end{equation}

\begin{remark}[Special case $k=2$]
    If $k=2$, there is a shortcut: A word does not contain a soft streak of length $2$ if and only if it is strictly decreasing. The number of strictly decreasing words of length $s$ is simply the number of combinations of $\{1,2,\ldots,n\}$ of length $s$, i.e. $\binom{n}{s}$. Thus the generating function is given by:
    \begin{equation*}
        f(z)
        = \sum_{s=0}^\infty \binom{n}{s}z^s
        = (1+z)^n
    \end{equation*}
    This agrees with \eqref{soft_streak:generating_function_eq} and \eqref{soft_streak:generating_function_alt_eq} as $\psi_{2,r}=(-1)^r$, $\omega_2=-1$ and $\Bs(n,2,r)=\binom{n}{r}$.\\
    
    Note that this series terminates at $s=n$. This makes sense, since if $s>n$, you cannot avoid a soft streak of length $2$ (the longest strictly decreasing word is $(n,n-1,\ldots,2,1)$ which has length $n$).
\end{remark}

\section{Application to Random Sampling}

We now return to the general Goulden-Jackson method. Suppose we sample letters randomly and uniformly from an alphabet of size $n$ (or equivalently, roll a fair $n$-sided die) and write them down successively until we encounter a forbidden word. Then the expected number of letters we must draw has a simple expression in terms of the generating function obtained by the Goulden-Jackson method.\\

Suppose $f(z)=\alpha_0+\alpha_1z+\alpha_2z^2+\cdots$ is the generating function for the number of words of any given length that do not contain any forbidden words. Also suppose $X$ is the random variable denoting the number of letters drawn until we encounter a forbidden word. Then for any $s\in\N_0$, $P(X>s)$ is the probability that we did \emph{not} encounter a forbidden word among the first $s$ letters. This yields:
\begin{equation}\label{expected_value_from_generating_function_eq}
    \E(X)
    = \sum_{s=0}^{\infty} P(X>s)
    = \sum_{s=0}^{\infty} \frac{\alpha_s}{n^s}
    = f\qty(\frac{1}{n})
\end{equation}
By \Cref{goulden_jackson:radius_of_convergence}, the power series of $f$ has radius of convergence greater than $\frac{1}{n}$, so it converges at $z=\frac{1}{n}$ and so \eqref{expected_value_from_generating_function_eq} is valid.

\begin{theorem}[Expected number of draws to get a streak]\label{streak:expected_value}
    Suppose consecutive letters are drawn from $\{1,2,\ldots,n\}$ until a streak of $k$ is obtained. Then the expected number of letters drawn is given by:
    \begin{align}\label{streak:expected_value_eq}
        E(n,k)
        &= \frac{1}{\Disp\sum_{r=0}^n \psi_{k,r}\binom{n}{r}n^{-r}} \\
    \label{streak:expected_value_eq2}
        &= \frac{k}{\Disp\sum_{s=1}^{k-1} \qty(1-\omega_k^{-s})\qty(1+\frac{\omega_k^s}{n})^n}
    \end{align}
\end{theorem}
These follow directly from \eqref{expected_value_from_generating_function_eq} by substituting $z=\frac{1}{n}$ into \eqref{streak:generating_function_eq} and \eqref{streak:generating_function_alt_eq} respectively.\\

We also have the following conjecture as a corollary of \eqref{expected_value_from_generating_function_eq}:

\begin{conjecture}[Conjecture \normalfont\sffamily(Expected number of draws to get a soft streak)]\label{soft_streak:expected_value}
    Suppose consecutive letters are drawn from $\{1,2,\ldots,n\}$ until a soft streak of $k$ is obtained. Then the expected number of letters drawn is given by:
    \begin{align}\label{soft_streak:expected_value_eq}
        E_{\textsf{soft}}(n,k)
        &= \frac{\qty(1-n^{-k})^n}{\Disp\sum_{r=0}^{(k-1)n} \psi_{k,r}\Bs(n,k,r)n^{-r}} \\
    \label{soft_streak:expected_value_eq2}
        &= \frac{k}{\Disp\sum_{s=1}^{k-1} \qty(1-\omega_k^{-s})\qty(1-\frac{\omega_k^s}{n})^{-n}}
    \end{align}
\end{conjecture}
These follow directly from \eqref{expected_value_from_generating_function_eq} by substituting $z=\frac{1}{n}$ into \eqref{soft_streak:generating_function_eq} and \eqref{soft_streak:generating_function_alt_eq} respectively.

\section{Continuous limit}

We now investigate the limit as $n\to\infty$ (recall that $n$ is the size of alphabet, or equivalently, the number of sides of the die). Intuitively, this should be equivalent to replacing the die with a spinner that generates random numbers in a continuous distribution on an interval. The exact distribution is not important in the case of a streak (as long as it is absolutely continuous), as we only care about monotonicity of the values we draw and not the values themselves.\\

Taking the limit as $n\to\infty$ in both \eqref{streak:expected_value_eq2} and \eqref{soft_streak:expected_value_eq2} yield the same result:
\begin{align*}
    \limn E(n,k) = \limn E_{\textsf{soft}}(n,k) = \frac{k}{\Disp\sum_{s=1}^{k-1} e^{\omega_k^s}\qty(1-\omega_k^{-s})}
\end{align*}
The case $k=3$ also appears in \cite[\S1.2]{lin}.\\

We now prove that this is indeed the result of the continuous case.\\

We wish to find $\mu_k$, the expected number of draws from an absolutely continuous distribution (which we will take to be uniform on $[0,1]$ for simplicity) to get a streak of $k$.\\

For each $r\in\{1,2,\ldots,k-1\}$, define $y_r(x)$ as the expected number of draws to get a streak of $k$, given that you currently have a streak of $r$ ending with $x$. This means that the last $r$ numbers you drew were in increasing order (but not the last $r+1$ numbers\footnotetext{This includes the case where you have only drawn $r$ numbers so far.}) and the last number you drew was $x$. For example, if you have drawn $0.43, 0.35, 0.67, 0.84$, then you currently have a streak of $3$ ending with $0.84$.

\begin{theorem}
    For all $k\in\N$, $k\ge 2$, we have:
    \begin{equation*}
        \mu_k = \frac{k}{\Disp\sum_{s=1}^{k-1} e^{\omega_k^s}\qty(1-\omega_k^{-s})}
    \end{equation*}
\end{theorem}
\begin{proof}
    Suppose you currently have a streak of $r\in\{1,2,\ldots,k-1\}$ ending with $x\in[0,1]$.
    \begin{itemize}
        \item If your next draw is $t<x$ (this has probability $x$), you have to start over (i.e. you have a streak of $1$ ending with $t$). Thus you expect to make $y_1(t)$ more draws.
        \item If your next draw is $t>x$ (this has probability $1-x$), you now have a streak of $r+1$ ending with $t$. Thus you expect to make $y_{r+1}(t)$ more draws.
    \end{itemize}
    And of course, the first number we draw gives us a streak of $1$ ending with that number. This yields the following system of integral equations:
    \begin{align}\label{continuous_streak:integral_equations}
        \begin{split}
            \mu_k &= 1 + \int_0^1 y_1(t)\dt \\
            y_1(x) &= 1 + \int_0^x y_1(t)\dt + \int_x^1 y_2(t)\dt \\
            y_2(x) &= 1 + \int_0^x y_1(t)\dt + \int_x^1 y_3(t)\dt \\
            & \vdotswithin{=} \\
            y_{k-2}(x) &= 1 + \int_0^x y_1(t)\dt + \int_x^1 y_{k-1}(t)\dt \\
            y_{k-1}(x) &= 1 + \int_0^x y_1(t)\dt
        \end{split}
    \end{align}
    Differentiating, we get:
    \begin{align}\label{continuous_streak:differential_equations}
        \begin{split}
            y_1'(x) &= y_1(x) - y_2(x) \\
            y_2'(x) &= y_1(x) - y_3(x) \\
            & \vdotswithin{=} \\
            y_{k-2}'(x) &= y_1(x) - y_{k-1}(x) \\
            y_{k-1}'(x) &= y_1(x)
        \end{split}
    \end{align}
    This is a system of $k-1$ linear ODEs in $y_1,y_2,\ldots,y_{k-1}$. As for the initial/boundary conditions, substituting $x=1$ into \eqref{continuous_streak:integral_equations} yields:
    \begin{equation}\label{continuous_streak:initial_conditions}
        y_1(1) = y_2(1) = \cdots = y_{k-1}(1) = \mu_k
    \end{equation}
    We must also have $y_{k-1}(0)=1$, since if you currently have a streak of $k-1$ ending with $0$, you are guaranteed to get a streak of $k$ with the next draw\footnotemark.\\
    
    If we ``normalize'' the functions $y_1,y_2,\ldots,y_{k-1}$ by setting $u_r=\frac{1}{\mu_k}y_r$, we get a purely initial value problem:
    \begin{align}\label{continuous_streak:differential_equations_normalized}
        \begin{aligned}
            u_1'(x) &= u_1(x) - u_2(x) \\
            u_2'(x) &= u_1(x) - u_3(x) \\
            & \vdotswithin{=} \\
            u_{k-2}'(x) &= u_1(x) - u_{k-1}(x) \\
            u_{k-1}'(x) &= u_1(x)
        \end{aligned}
        &&
        u_1(1) = u_2(1) = \cdots = u_{k-1}(1) = 1
    \end{align}
    And $\mu_k=\frac{1}{u_{k-1}(0)}$.\\
    
    The solution of \eqref{continuous_streak:differential_equations_normalized} is given by:
    \begin{equation}
        u_r(x) = \frac{1}{k}\sum_{s=1}^{k-1} \qty(1-\omega_k^{rs})e^{\omega_k^s(1-x)}
        \tag{$r=1,2,\ldots,k-1$}
    \end{equation}
    Substituting $r=k-1$ and $x=0$ and inverting yields:
    \begin{equation*}
        \mu_k = \frac{1}{u_{k-1}(0)} = \frac{k}{\Disp\sum_{s=1}^{k-1} e^{\omega_k^s}\qty(1-\omega_k^{-s})} \qedhere
    \end{equation*}
\end{proof}
\footnotetext{Ok, this is not \emph{technically} correct, but since the probability of drawing the same number more than once is zero, we can safely ignore this issue.}

\subsubsection*{Explicit formulas for small values of $k$}

For $k=2$:
\begin{align*}
    u_1(x) = e^{x-1} && \mu_2 = e \approx 2.718282
\end{align*}

For $k=3$:
\begin{align*}
    y_1(x) = \frac{e^{\tfrac{x}{2}}\cos(\frac{\sqrt{3}}{2}(x-1)+\frac{\pi}{6})}{\cos(\frac{\sqrt{3}}{2}+\frac{\pi}{6})}
    &&
    y_2(x) = \frac{e^{\tfrac{x}{2}}\cos(\frac{\sqrt{3}}{2}(x-1)-\frac{\pi}{6})}{\cos(\frac{\sqrt{3}}{2}+\frac{\pi}{6})}
    &&
    \mu_3
    = \frac{\sqrt{3e}}{2\cos(\dfrac{\sqrt{3}}{2}+\dfrac{\pi}{6})}
    \approx 7.924372
\end{align*}
For $k=4$:
\begin{align*}
    y_1(x) &= \frac{e^x-e\sin(x-1)+e\cos(x-1)}{1-e(\sin(1)-\cos(1))} \\
    y_2(x) &= \frac{2e\cos(x-1)}{1-e(\sin(1)-\cos(1))} \\
    y_3(x) &= \frac{e^x+e\sin(x-1)+e\cos(x-1)}{1-e(\sin(1)-\cos(1))} \\
    \mu_4 &= \frac{2e}{1-e(\sin(1)-\cos(1))}
    \approx 29.980170
\end{align*}

For higher values of $k$, we have:

\begin{center}\ttfamily
    \newcommand*{\highdigits}[1]{\textcolor{blue!50!black}{#1}}
    \begin{tabular}{|r|r|}
        \multicolumn{1}{c}{$k$} & \multicolumn{1}{c}{$\mu_k$ (to $60$ decimal places)} \\\hline
        2 & \highdigits{2.}718281828459045235360287471352662497757247093699959574966968 \\\hline
        3 & \highdigits{7.9}24372434513184628799810694208415749584787751007803272302703 \\\hline
        4 & \highdigits{29.9}80170111893315322399955239983927299141905015329020581286930 \\\hline
        5 & \highdigits{143.99}4805367336197413565390684410970568671218564664180561091334 \\\hline
        6 & \highdigits{839.99}8640250736936484708350540147060715630641050163631672497854 \\\hline
        7 & \highdigits{5759.999}644799666083821542490440199192540829862157792061849437504 \\\hline
        8 & \highdigits{45359.9999}07445495866458458747887235726727998085580263917214428834 \\\hline
        9 & \highdigits{403199.9999}75944248391489643004925481397778149807848098524214919893 \\\hline
        10 & \highdigits{3991679.99999}3762687261148961681589183499555116134021765003499924135 \\\hline
        11 & \highdigits{43545599.99999}8386322252253123523422442886632124901510775067769335874 \\\hline
        12 & \highdigits{518918399.999999}583357370408117814222274969576150520059702617407517123 \\\hline
        13 & \highdigits{6706022399.999999}892620702958568083451300164881681275187541538503077141 \\\hline
        14 & \highdigits{93405311999.9999999}72371178409026358124994762865268191373551870561927755 \\\hline
        15 & \highdigits{1394852659199.99999999}2901670495910427777148139210921443892403824212534308 \\\hline
        16 & \highdigits{22230464255999.99999999}8178784422334973035853644616759173438447029466326330 \\\hline
        17 & \highdigits{376610217983999.999999999}533307881630474553527658533401701129555158955094942 \\\hline
        18 & \highdigits{6758061133823999.999999999}880543479711637568193161366429059063487770379745148 \\\hline
        19 & \highdigits{128047474114559999.9999999999}69454970307684867213252426425625170841619171498343 \\\hline
        20 & \highdigits{2554547108585471999.99999999999}2197048664880838487190307346342867240167220649762 \\\hline
        21 & \highdigits{53523844179886079999.99999999999}8008422482625283883246379136011471180505449109885 \\\hline
        22 & \highdigits{1175091669949317119999.999999999999}492092868675029353971166800320262160032013123068 \\\hline
        23 & \highdigits{26976017466662584319999.999999999999}870566722275066805564948431276562750692499255348 \\\hline
        24 & \highdigits{646300418472124415999999.9999999999999}67038630823130763357357810721382325899452430946 \\\hline
        25 & \highdigits{16131658445064225423359999.99999999999999}1611525434893702034663107782646382044952491906 \\\hline
        26 & \highdigits{418802671169936621567999999.99999999999999}7866473861936721552720430101299686177977844790 \\\hline
        27 & \highdigits{11292160911544957796351999999.999999999999999}457665584258861854728186147710347961739620511 \\\hline
        28 & \highdigits{315777214062132212662271999999.999999999999999}862213708621620749813293533566716693152765422 \\\hline
        29 & \highdigits{9146650338351415815045119999999.9999999999999999}65011240223148822984088726555162948340854195 \\\hline
        30 & \highdigits{274094621805930760590852095999999.99999999999999999}1119291928769895698939605832111440868044313 \\\hline
    \end{tabular}
\end{center}

This suggests that $\mu_k=k!+(k-1)!+o(1)$ as $k\to\infty$. We were not able to prove this claim.

\appendix

\section{Lower bound for the radius of convergence}\label{radius_of_convergence:lower_bound}

In this section:
\begin{itemize}
    \item The alphabet will always be $\As=\{1,2,\ldots,n\}$, where $n\in\N$ is fixed.
    \item $f_\Fs$ denotes the generating function corresponding to a set $\Fs$ of forbidden words with letters in $\As$.
    \item $\Radius(f)$ denotes the radius of convergence of a function $f$ expressed as a power series, which will always be centered at $0$.
\end{itemize}

\vspace{10pt}

If $\Fs$ is empty, we have $f_\Fs(z)=\frac{1}{1-nz}$, which has radius of convergence $\frac{1}{n}$. The goal of this section is to prove that if $\Fs$ is non-empty, the radius of convergence is \emph{strictly greater} than $\frac{1}{n}$.

\begin{lemma}\label{radius_of_convergence:inequality}
    Suppose $A(z)=\sumn[0] \alpha_nz^n$ and $B(z)=\sumn[0] \beta_nz^n$ are two power series such that $0\le\alpha_n\le\beta_n$ for all $n\in\N_0$. Then the radius of convergence of $A(z)$ is at least as large as that of $B(z)$.
\end{lemma}
\begin{proof}
    \begin{equation*}
    0\le\alpha_n\le\beta_n \implies 0\le\limsupn\alpha_n\le\limsupn\beta_n \implies \frac{1}{\limsupn\alpha_n} \ge \frac{1}{\limsupn\beta_n} \qedhere
    \end{equation*}
\end{proof}

Which set of forbidden words is the least ``forbidding''? In other words, which $\Fs$ gives rise to the generating function with the largest possible coefficients? Of course, this question is ill-defined as we can always increase the lengths of the forbidden words. However, if we require the longest forbidden word to have some fixed length $k$, then it is clear the answer involves a single forbidden word of length $k$ (since we can also increase the coefficients by removing forbidden words). We also want this word to front-run itself as much as possible, to maximize the number of ways to avoid it as a subword. This leads to the following definition:

\begin{definition}
    Suppose $k\ge 2$. The \textbf{minimal forbidden set} of \textbf{size} $k$ is given by $\Ms(k)=\{(1,1,\ldots,1)\}$.
\end{definition}
\begin{remark}
    The choice of $1$ is not important, what matters is that there is exactly one forbidden word and this word comprises the same letter repeated $k$ times.
\end{remark}

\begin{proposition}\label{minimal_forbidden_set:generating_function}
    Suppose $k\in\N$, $k\ge 2$. Then the generating function corresponding to $\Ms(k)$ is given by:
    \begin{equation}\label{minimal_forbidden_set:generating_function_eq}
        f_{\Ms(k)}(z) = \frac{1+z+z^2+\cdots+z^{k-1}}{1-(n-1)\qty(z+z^2+\cdots+z^k)}
    \end{equation}
\end{proposition}
\begin{proof}
    Since the only forbidden word is $(1,1,\ldots,1)$, which front-runs itself with overlaps of $1,2,\ldots,k-1$, we have:
    \begin{align*}
        W(1,1,\ldots,1)
        &= -z^k - z^{k-1}W(1,1,\ldots,1) - \cdots - z^2W(1,1,\ldots,1) - zW(1,1,\ldots,1) \\
        &= -z^k - (z+z^2+\cdots+z^{k-1})W(1,1,\ldots,1) \\
        (1+z+z^2+\cdots+z^{k-1})W(1,1,\ldots,1)
        &= -z^k \\
        W(\Fs) = W(1,1,\ldots,1)
        &= -\frac{z^k}{1+z+z^2+\cdots+z^{k-1}}
    \end{align*}
    This yields:
    \begin{align*}
        f_{\Ms(k)}(z)
        &= \frac{1}{1-nz-W(\Fs)} \\
        &= \frac{1}{1-nz+\dfrac{z^k}{1+z+z^2+\cdots+z^{k-1}}} \\
        &= \frac{1+z+z^2+\cdots+z^{k-1}}{(1-nz)\qty(1+z+z^2+\cdots+z^{k-1})+z^k} \\
        &= \frac{1+z+z^2+\cdots+z^{k-1}}{1-(n-1)\qty(z+z^2+\cdots+z^k)} \qedhere
    \end{align*}
\end{proof}

\begin{lemma}\label{minimal_forbidden_set:nonvanishing}
    Suppose $n,k\in\N$, $k\ge 2$. Then $\widehat{M}_{n,k}(z)>0$ for all $z\in\qty[0,\frac{1}{n}]$.
\end{lemma}
\begin{proof}
    Define $p(z)=1+z+z^2+\cdots+z^{k-1}$ and $q(z)=1-(n-1)\qty(z+z^2+\cdots+z^k)$ as the numerator and denominator respectively of \eqref{minimal_forbidden_set:generating_function_eq}. Clearly $p(z)>0$ for all $z\ge 0$. We also have:
    \begin{equation*}
        q'(z) = -(n-1)\qty(1+2z+3z^2+\cdots+kz^{k-1}) < 0 \text{ for all } z \ge 0
    \end{equation*}
    Thus $q$ is strictly decreasing on $[0,\infty)$. Furthermore:
    \begin{equation*}
        q\qty(\frac{1}{n})
        = 1-(n-1)\qty(\frac{1}{n}+\frac{1}{n^2}+\cdots+\frac{1}{n^k})
        = 1-(n-1)\qty(\frac{1-n^{-k}}{n-1})
        = n^{-k} > 0
    \end{equation*}
    Thus $q(z)>0$ for all $z\in\qty[0,\frac{1}{n}]$, and so $\widehat{M}_{n,k}(z)>0$ for all $z\in\qty[0,\frac{1}{n}]$.
\end{proof}

\begin{theorem}[Vivanti-Pringsheim Theorem]\label{vivanti_pringsheim}
    Suppose $f(z)=\sum_{s=0}^\infty \alpha_s z^s$ is a power series such that $\alpha_s\ge 0$ for all $s\in\N_0$, and $R$ is the radius of convergence of $f(z)$. Then $R$ is a singularity of $f(z)$.
\end{theorem}
See \cite[Theorem IV.6, Page 240]{flajolet} or \cite[Page 235]{remmert} for a proof.

\begin{theorem}[Strict lower bound on the radius of convergence]\label{goulden_jackson:radius_of_convergence}
    Suppose $\Fs$ is a non-empty set of words with letters in $\As$. Then $\Radius(f_\Fs)>\frac{1}{n}$.
\end{theorem}
\begin{proof}
    Suppose $k$ is the length of the longest forbidden word, and $f_\Fs(z)=\sum_{s=0}^\infty \alpha_sz^s$ and $f_{\Ms(k)}(z)=\sum_{s=0}^\infty \beta_sz^s$. Then we have $0\le\alpha_s\le\beta_s$ for all $s\in\N_0$, so by \Cref{radius_of_convergence:inequality}, we have $\Radius(f_\Fs)\ge\Radius(f_{\Ms(k)})$. By the \hyperref[vivanti_pringsheim]{Vivanti-Pringsheim Theorem}, $\Radius(f_{\Ms(k)})$ is a singularity of $f_{\Ms(k)}$, and so by \Cref{minimal_forbidden_set:nonvanishing}, we have $\Radius(f_{\Ms(k)})>\frac{1}{n}$. Thus $\Radius(f_\Fs)>\frac{1}{n}$.
\end{proof}
\begin{remark}
    
\end{remark}

\printbibliography

\end{document}